\documentclass[11pt,final]{amsart}
\usepackage[utf8]{inputenc}
\usepackage[T1]{fontenc}
\usepackage[letterpaper,centering]{geometry}
\usepackage{lmodern}
\usepackage{eucal}
\usepackage{mathrsfs}
\usepackage{enumerate}
\usepackage{amsmath,amssymb,amsfonts}
\usepackage{xcolor}
\definecolor{darkblue}{rgb}{0,0,0.3}
\definecolor{darkgreen}{rgb}{0,0.4,0}
\usepackage[colorlinks=true,linkcolor=darkblue, citecolor=darkgreen, unicode,pdfborder={0 0 0}]{hyperref}
\usepackage[ps,arrow,matrix,tips,line, curve]{xy}
{\setbox0\hbox{$ $}}\fontdimen16\textfont2=\fontdimen17\textfont2
\entrymodifiers={+!!<0pt,\the\fontdimen22\textfont2>}
\SelectTips{cm}{11}
\usepackage{enumitem}
\setlist[enumerate]{label={\upshape(\arabic*)},topsep=.7ex, leftmargin=*}
\setlist[itemize]{leftmargin=*}

\usepackage{ifdraft}
\ifdraft{
\usepackage[notcite,notref]{showkeys}
\addtolength{\hoffset}{1.5cm}

}{}

\theoremstyle{plain}

\newtheorem{thmintro}{Theorem}

\newtheorem{thm}{Theorem}[section]

\theoremstyle{definition}
\newtheorem{question}[thm]{Question}
\newtheorem{questions}[thm]{Questions}

\newtheorem{defn}[thm]{Definition}

\numberwithin{equation}{section}

\input cyracc.def
\DeclareFontFamily{U}{russian}{}
\DeclareFontShape{U}{russian}{m}{n}
        { <5><6> wncyr5
        <7><8><9> wncyr7
        <10><10.95><12><14.4><17.28><20.74><24.88> wncyr10 }{}
\DeclareSymbolFont{Russian}{U}{russian}{m}{n}
\DeclareSymbolFontAlphabet{\mathcyr}{Russian}
\makeatletter
\let\@math@cyr\mathcyr
\renewcommand{\mathcyr}[1]{\@math@cyr{\cyracc #1}}
\makeatother
\newcommand{\ci}{\sC^\infty}

\newcommand{\rec}{\text{rec}}

\makeatletter
\def\myrightarrow{{\setbox\z@\hbox{$\rightarrow$}\dimen0\ht\z@\multiply\dimen0 6\divide\dimen0 10\ht\z@\dimen0\box\z@}}
\def\myrightarrowfill@{\arrowfill@\relbar\relbar\myrightarrow}
\def\myleftarrow{{\setbox\z@\hbox{$\leftarrow$}\dimen0\ht\z@\multiply\dimen0 6\divide\dimen0 10\ht\z@\dimen0\box\z@}}
\def\myleftarrowfill@{\arrowfill@\myleftarrow\relbar\relbar}
\newcommand{\myxrightarrow}[2][]{\ext@arrow 0359\myrightarrowfill@{#1}{#2}}
\newcommand{\myxleftarrow}[2][]{\ext@arrow 3095\myleftarrowfill@{#1}{#2}}
\makeatother
\newcommand{\mtilde}{{\mathchoice
    {\widetilde{m}}
    {\widetilde{m}}
    {\rlap{$\scriptscriptstyle{m}$}\vphantom{\raise0pt\hbox{$m$}}\smash{\lower2.5pt\hbox{$\scriptscriptstyle\widetilde{\phantom{\scriptscriptstyle{m}}}$}}}
    {\rlap{$\scriptscriptstyle{m}$}\vphantom{\raise.2pt\hbox{$m$}}\smash{\lower2.05pt\hbox{$\scriptscriptstyle\widetilde{\phantom{\scriptscriptstyle{m}}}$}}}}}

\newcommand{\Mtilde}{{\mathchoice
    {\rlap{$M$}\mkern1mu\smash[b]{\lower.5pt\hbox{$\widetilde{\phantom{M}}$}}\mkern-1mu}
    {\rlap{$M$}\mkern1mu\smash[b]{\lower.5pt\hbox{$\widetilde{\phantom{M}}$}}\mkern-1mu}
    {\rlap{$\scriptstyle{M}$}\mkern1mu\smash[b]{\lower.5pt\hbox{$\widetilde{\phantom{\scriptstyle{M}}}$}}\mkern-1mu}
    {\widetilde{M}}}}

\newcommand{\et}{{\text{ét}}}
\newcommand{\sC}{{\mathscr C}}

\newcommand{\sB}{{\mathscr B}}

\newcommand{\sX}{{\mathscr X}}
\newcommand{\sZ}{{\mathscr Z}}
\newcommand{\sY}{{\mathscr Y}}
\newcommand{\A}{{\mathbf A}}

\renewcommand{\C}{{\mathbf C}}

\newcommand{\cF}{\mathrm F}

\newcommand{\cFplus}{\mathrm F_+}
\newcommand{\cFconst}{\mathrm F_\const}

\renewcommand{\P}{{\mathbf P}}
\newcommand{\Q}{{\mathbf Q}}
\newcommand{\R}{{\mathbf R}}
\renewcommand{\S}{{\mathbf S}}

\newcommand{\Z}{{\mathbf Z}}

\newcommand{\Qp}{{\Q_{p}}}
\newcommand{\CH}{\mathrm{CH}}

\newcommand{\nr}{\mathrm{nr}}
\newcommand{\alg}{\mathrm{alg}}
\newcommand{\const}{\mathrm{const}}

\newcommand{\Gm}{\mathbf{G}_\mathrm{m}}

\newcommand{\Gal}{\mathrm{Gal}}
\newcommand{\SL}{\mathrm{SL}}
\newcommand{\Pic}{\mathrm{Pic}}
\newcommand{\Br}{\mathrm{Br}}

\renewcommand{\phi}{\varphi}
\renewcommand{\emptyset}{\varnothing}

\newcommand{\chapeau}{{\rlap{\smash{\hbox{\lower4pt\hbox{\hskip1pt$\widehat{\phantom{u}}$}}}}}}
\newcommand{\Picplushat}{\Pic_+^{{\smash{\hbox{\lower4pt\hbox{\hskip0.4pt$\widehat{\phantom{u}}$}}}}}}
\newcommand{\PicplusAhat}{\Pic_{+,\A}^{{\smash{\hbox{\lower4pt\hbox{\hskip.4pt$\widehat{\phantom{u}}$}}}}}}

\newcommand{\Pichat}{\Pic^{{\smash{\hbox{\lower4pt\hbox{\hskip0.4pt$\widehat{\phantom{u}}$}}}}}}

 \makeatletter
 \renewcommand{\tocsection}[3]{%
   \indentlabel{\@ifnotempty{#2}{\bfseries\ignorespaces#1 #2\quad}}\bfseries#3}

 \renewcommand{\tocsubsection}[3]{%
   \indentlabel{\@ifnotempty{#2}{\hspace{1.6em}\ignorespaces#1 #2\quad}}#3}
 \makeatother

\advance\topmargin -1.5mm
\advance\textheight 3mm

\makeatletter\let\@wraptoccontribs\wraptoccontribs\makeatother

\date{October 29th, 2021}
\title{Some aspects of rational points and rational curves}

\author{Olivier Wittenberg}
\address{Institut Galil\'ee, Universit\'e Sorbonne Paris Nord, 99~avenue Jean-Baptiste Cl\'ement, 93430 Villetaneuse, France}
\email{wittenberg@math.univ-paris13.fr}

\begin{document}
\begin{abstract}
Various methods have been used to construct rational points and rational
curves on rationally connected algebraic varieties.  We survey recent advances in two of them,
the descent and the fibration method, in a number-theoretical context
(rational points over number fields) and in an algebro-geometric one (rational curves on real varieties),
and discuss the question of rational points over function fields of $p$\nobreakdash-adic curves.
\end{abstract}

\maketitle
\section{Introduction}

Let~$X$ denote an algebraic variety over a field~$k$ and~$X(k)$
the set of its rational points.

The search for explicit descriptions of the set $X(k)$
 when~$k$ is a number field
is one of the oldest themes of number theory.
A modern point of view on this problem
consists in embedding~$X(k)$ diagonally
into the topological space~$X(\A_k)$ of adelic points of~$X$
and attempting to identify its topological closure.
By general principles that were formulated by Lang
after the works of Mordell, Weil and Siegel,
the answer is expected
 to  depend in a crucial manner on the geometry of~$X$.
For instance, assuming that~$X$ is smooth and projective and that an embedding $k \hookrightarrow \C$ is given,
the set~$X(k)$ is conjectured to be finite if the complex variety $X_{\C}$
is \emph{hyperbolic} (see \cite{langhigherdim}). One may then seek to count, list or bound its elements.
At the other end of the spectrum, if~$X_{\C}$ is a \emph{rationally connected} smooth projective
variety in the sense
of Campana \cite{campana} and Kollár--Miyaoka--Mori \cite{kmm}, then
one expects that the set~$X(k)$ is Zariski dense in~$X$ whenever it is nonempty; more precisely,
by
a conjecture of Colliot-Thélène,
the closure of~$X(k)$ in $X(\A_k)$ should coincide in this case with the Brauer--Manin set $X(\A_k)^{\Br(X)}$
defined by Manin~\cite{maninicm}.
This far-reaching conjecture encompasses in particular the inverse Galois problem, and its refinement
the Grunwald problem (see \cite{ekedahl}, \cite[\textsection3.5]{serretopics}, \cite{harariquelques},
\cite{dlan}).

Criteria for the existence of rational points on~$X$ are also of relevance outside of number theory,
when~$k$ is no longer assumed to be a number field.
For instance, the Graber--Harris--Starr theorem \cite{ghs}, a central result in the theory of
rational curves on complex algebraic varieties, is equivalent to the statement that $X(k)\neq\emptyset$
if $k$ is the function field of a complex
curve
and~$X$ is a rationally connected variety.
(We say that~$X$ is \emph{rationally connected}
to mean that
for any algebraically closed field extension~$K$ of~$k$,
the variety~$X_K$ over~$K$ is rationally connected in the sense of \cite{campana,kmm}.)
As another example, if~$X$
is a real algebraic variety with no real point
and~$k$ denotes
 the function field of the real conic given by $x^2+y^2=-1$,
the existence of a geometrically rational curve on~$X$---a property conjectured  by Kollár
to hold
whenever~$X$ is
a positive-dimensional
rationally connected variety---is equivalent to the statement
that $X(k)\neq\emptyset$.

The results we discuss in this expository article concern the existence of rational points
in two very distinct contexts, leading
 to the following two concrete theorems.
As we shall see, their proofs roughly follow, perhaps somewhat surprisingly, a common general strategy.

\newcommand{\citethmintroA}{\cite{hwzceh}}
\begin{thmintro}[see \citethmintroA]
\label{thmintro:A}
Let~$G$ be a finite nilpotent group.
Let~$k$ be a number field.
\begin{enumerate}
\item\label{item1} There exist Galois extensions~$K/k$ with Galois group~$G$.
\item 
If $v_1,\dots,v_n$ are pairwise distinct places of~$k$ none of which is a finite place dividing the order of~$G$,
and $w_1,\dots,w_n$ are places of~$K$ above $v_1,\dots,v_n$,
then in~\ref{item1}, one can require that the extensions $K_{w_i}/k_{v_i}$ be isomorphic to any prescribed
collection of
Galois extensions of~$k_{v_1},\dots,k_{v_n}$ whose Galois groups are subgroups of~$G$.
\end{enumerate}
\end{thmintro}

\newcommand{\citethmintroB}{\cite{bwtight}}
\begin{thmintro}[see \citethmintroB]
\label{thmintro:B}
Let~$X$ be a smooth, proper variety over~$\R$.  Let $\varepsilon:\S^1 \to X(\R)$ be a continuous map.
Assume that~$X$ is birationally equivalent to a homogeneous space of a linear algebraic group over~$\R$.
Then there exist morphisms of algebraic varieties $\P^1_\R \to X$
that induce maps $\P^1(\R)=\S^1 \to X(\R)$ arbitrarily close to~$\varepsilon$ in the compact-open
topology.
\end{thmintro}

Theorem~\ref{thmintro:A}~(1) was first proved by Shafarevich in his seminal work on the inverse Galois problem
for solvable groups (see \cite[Chapter~IX, \textsection6]{nsw};
it should be noted that
nilpotent groups form the most difficult case in his proof);
the proof given in~\cite{hwzceh} is independent from his and has a geometric flavour.
Theorem~\ref{thmintro:A}~(2), on the other hand,
was new in~\cite{hwzceh} and was not accessible with Shafarevich's methods.

As far as we know,
Theorem~\ref{thmintro:B} might hold under the sole assumption that~$X$ is rationally connected.
This is a question we raise in \cite{bwtight}.
Theorem~\ref{thmintro:B} provides the first examples of a positive answer to it for varieties
that are not
$\R$\nobreakdash-rational (indeed, not even $\C$\nobreakdash-rational).
For $\R$\nobreakdash-rational varieties, the conclusion of Theorem~\ref{thmintro:B} was previously
shown, by Bochnak and Kucharz~\cite{BKrat}, to follow from the Stone--Weierstrass theorem.

The first step in the proofs of Theorem~\ref{thmintro:A} and Theorem~\ref{thmintro:B}
consists in strengthening and reformulating
the desired conclusion in terms of the existence of suitable rational points on suitable varieties over suitable fields.
In the case of Theorem~\ref{thmintro:A},
the varieties in question are homogeneous spaces of~$\SL_n$ over number fields; for the proof, though
not for the statement, it is crucial to \emph{not} restrict to homogeneous spaces that have rational points
(i.e.\ to homogeneous spaces of the form $\SL_n/G$).
In the case of Theorem~\ref{thmintro:B}, the varieties in question are homogeneous spaces of linear algebraic
groups, over the rational function field~$\R(t)$; for the proof, though not for the statement,
it is crucial to \emph{not} restrict to
homogeneous spaces or algebraic groups that are defined over~$\R$.
In the remainder of the proofs of Theorems~\ref{thmintro:A} and~\ref{thmintro:B}, one establishes
the validity of these strengthened formulations by combining geometric \emph{dévissages} of the underlying
algebraic varieties
with two general tools: the descent method and the fibration method.
The fibration method, whose first instance can be found in the work of Hasse on the local-global principle
for quadratic forms, consists
in reducing
the desired property for a variety~$V$
 endowed with a morphism
$p:V \to B$ with geometrically irreducible generic fibre
to the same property for~$B$ and for a collection of smooth fibres of~$p$.
The descent method, which goes back to Fermat, attempts to reduce the desired property for a variety~$V$
 endowed with a torsor $p:W\to V$ under a (possibly disconnected)
linear algebraic group over~$k$
to the same property for~$W$ and for all of its twists.
It was developed in the context of elliptic curves, for torsors under finite abelian groups,
by Mordell, Cassels and Tate, and the set-up was later extended to torsors under positive-dimensional
linear algebraic groups by Colliot-Thélène and Sansuc, Skorobogatov, Harari.

We take Theorems~\ref{thmintro:A} and~\ref{thmintro:B} as excuses
leading us to the general study of rational points on rationally
connected varieties defined over number fields or over function fields of real curves.
We discuss recent advances in the fibration and descent methods in these two contexts
in~\textsection\ref{sec:numberfields} and in~\textsection\ref{sec:tight},
stating along the way the main open questions that surround  Theorems~\ref{thmintro:A} and~\ref{thmintro:B}
and their proofs. We then turn,
in~\textsection\ref{sec:otherfields},
to function fields of $p$\nobreakdash-adic curves,
and speculate about the existence of a $p$\nobreakdash-adic analogue
of the ``tight approximation'' property discussed in~\textsection\ref{sec:tight}
that would enable one to exploit fibration and descent methods
in the study of rational curves over $p$\nobreakdash-adic fields
and more generally of rational points over function fields of $p$\nobreakdash-adic curves.

\bigskip
\emph{Acknowledgements.}
I am grateful to
Jean-Louis Colliot-Thélène,
Antoine Ducros,
Hélène Esnault,
János Kollár and
Ján Mináč for their comments on a first version of the text,
and to Olivier Benoist and to Yonatan Harpaz for the pleasant collaborations that
have led to
 the results reported on in this article.

\section{Solvable groups and the Grunwald problem in inverse Galois theory}
\label{sec:numberfields}

\subsection{Homogeneous spaces}

It is the following general theorem about the arithmetic of homogeneous spaces of linear algebraic groups
that underlies Theorem~\ref{thmintro:A}.

\begin{thm}
\label{thm:arithhomspace}
Let~$V$ be a homogeneous space of a connected
linear algebraic group~$L$ over a number field~$k$.
  Let~$X$ be a smooth compactification of~$V$.
Let $\bar v \in V(\bar k)$.  
Assume that the group of connected components~$G$ of the stabiliser of~$\bar v$
is supersolvable, in the sense that it possesses a normal series $1 = G_0 \lhd \dots \lhd G_m = G$ such
that the quotients~$G_{i+1}/G_i$ are cyclic while the subgroups~$G_i$ are normal in~$G$ and are stable under the
natural outer action of $\Gal(\bar k/k)$ on~$G$.
Then the subset~$X(k)$ is dense in $X(\A_k)^{\Br(X)}$.
\end{thm}

Here and elsewhere, by ``compactification of~$V$'', we mean a proper
variety over~$k$ that contains~$V$ as a dense open subset; we do not
require that the algebraic group~$L$ act on the compactification.
Examples of supersolvable groups with respect to
the trivial outer action of~$\Gal(\bar k/k)$ include finite nilpotent groups and dihedral
groups.  With a nontrivial outer action of~$\Gal(\bar k/k)$, however, even abelian groups need not be
supersolvable. Previous work of Borovoi~\cite{borovoi} nevertheless establishes
the conclusion of Theorem~\ref{thm:arithhomspace} in many cases
where the stabiliser
of~$\bar v$ is abelian but not necessarily supersolvable.

Theorem~\ref{thm:arithhomspace} can be found in
\cite[Théorème~B]{hwzceh}
in the particular case
where~$L$ is semi-simple simply connected and the stabiliser of~$\bar v$ is finite,
and in \cite[Corollary~4.5]{hwsupersolvable} in general.
To deduce Theorem~\ref{thmintro:A} from it, embed~$G$ into $\SL_n(k)$ for some~$n$,
take $L=\SL_n$ and $V=\SL_n/G$
and let~$H$ denote the set of points of~$V$ above which the fibre
of the étale cover $\pi:L \to V$ is irreducible.  The function field of the fibre of~$\pi$
above any rational point contained in~$H$ is a Galois extension of~$k$ with Galois group~$G$.
On the other hand, by a theorem of Ekedahl
\cite{ekedahl}, the density of $X(k)$
in $X(\A_k)^{\Br(X)}$ implies that of $X(k) \cap H$ in $X(\A_k)^{\Br(X)}$.
Thus, Theorem~\ref{thm:arithhomspace} ensures the existence of Galois extensions~$K/k$ with Galois
group~$G$ having a local behaviour prescribed by any element of the Brauer--Manin set $X(\A_k)^{\Br(X)}$;
that is, one may freely prescribe the completions of~$K$ at any finite set of places of~$k$,
as long as these prescriptions satisfy a certain global reciprocity
condition determined by~$\Br(X)$.
By a theorem of Lucchini Arteche
\cite[\textsection 6]{lucchiniunramifiedbrauer},
this reciprocity condition imposes, in fact,
no restriction at the places indicated in Theorem~\ref{thmintro:A}~(2).

\subsection{Geometry}
\label{subsec:geometry}

In the special case where $L=\SL_n$ and the stabiliser of~$\bar v$ is a finite group~$G$,
the geometry behind the proof of Theorem~\ref{thm:arithhomspace}
can be summarised with the following assertion: there exist an algebraic torus~$T$ over~$k$ and a
torsor $\bar Y \to X_{\bar k}$ under~$T_{\bar k}$
whose isomorphism class is invariant under $\Gal(\bar k/k)$,
such that for any torsor $Y \to X$ under~$T$ whose base change to~$X_{\bar k}$ is isomorphic to~$\bar Y$,
there exist a dense open subset $W \subseteq Y$ and
a smooth morphism
$p:W \to Q$ to a quasi-trivial torus~$Q$ (i.e.\ a torus of the form $R_{E/k}\Gm$ for a
nonzero étale $k$\nobreakdash-algebra~$E$) whose fibres
are homogeneous spaces of~$\SL_n$ with geometric stabiliser isomorphic to~$G_{m-1}$.
In addition, the morphism~$p$ admits a rational section over~$\bar k$.

This geometry is the key to a proof of Theorem~\ref{thm:arithhomspace} by
an induction on~$m$, at each step of which one applies
the descent method and the fibration method, in the form of
Theorem~\ref{th:descentnf}
and Theorem~\ref{th:fibrationnf}
below.
It should be noted that even if
$G$ is embedded into $\SL_n(k)$ and $V=\SL_n/G$, the homogeneous spaces of~$\SL_n$ that arise
as fibres of~$p$ need not possess rational points.  Thus, for the induction to be possible,
one cannot restrict
the
statement of Theorem~\ref{thm:arithhomspace} to homogeneous spaces of the form $\SL_n/G$,
even though only homogeneous spaces of this form
are relevant for Theorem~\ref{thmintro:A}.

\subsection{Descent}

The following theorem, which was established in \cite{hwzceh} and can also be deduced from
 \cite{caoapproxforte}, is the definitive statement of descent theory in the case of
smooth and proper rationally connected varieties over number fields.
For geometrically rational~$X$, this theorem is due to
Colliot-Thélène and Sansuc \cite{ctsandescent2}. The homogeneous spaces
of Theorem~\ref{thm:arithhomspace} are not geometrically rational in general
(Saltman, Bogomolov; see \cite{ctsmumbai}).

\begin{thm}
\label{th:descentnf}
Let~$X$ be a smooth and proper rationally connected variety over a number field~$k$.
Let~$T$ be a torus over~$k$ and $\bar Y \to X_{\bar k}$  a torsor under~$T_{\bar k}$
whose isomorphism class is invariant under $\Gal(\bar k/k)$.
Then
$$X(\A_k)^{\Br(X)}= \bigcup_{f:Y\to X}
f'\big(Y'(\A_k)^{\Br(Y')}\big)\rlap,$$
where the union ranges over the torsors $f:Y \to X$ under~$T$
whose base change to~$X_{\bar k}$ is isomorphic to~$\bar Y$,
and~$Y'$ denotes a smooth compactification of~$Y$ such that~$f$
extends to a morphism $f':Y' \to X$.
In particular, if $Y'(k)$ is dense in $Y'(\A_k)^{\Br(Y')}$ for every such~$f$,
then $X(k)$ is dense in $X(\A_k)^{\Br(X)}$.
\end{thm}

(To bridge the gap between Theorem~\ref{th:descentnf} and \cite[Théorème~2.1]{hwzceh},
one needs to know that $X(\A_k)^{\Br(X)}\neq\emptyset$ implies the existence
of at least one~$f$.
This goes back to~\cite{ctsandescent2} and follows from \cite[Theorem~3.3.1]{walb},
\cite[Proposition~2.2.5]{ctsandescent2}, \cite[(3.3)]{wittenbergslc}.)

\subsection{Fibration}

The following fibration theorem suffices for the proof of Theorem~\ref{thm:arithhomspace}.
It results from combining a descent
with the work of Harari \cite{harariduke} on the fibration method.

\begin{thm}
\label{th:fibrationnf}
Let $p:Z \to B$ be a dominant morphism between irreducible, smooth and proper varieties over a number field~$k$,
with rationally connected generic fibre.
Assume that
\begin{enumerate}
\item there exist dense open subsets $W \subset Z$ and $Q \subset B$ such that~$Q$ is a quasi-trivial torus
over~$k$ and~$p$ induces a smooth morphism $W \to Q$ with geometrically irreducible fibres;
\item the morphism $p$ admits a rational section over~$\bar k$;
\item
for all $b \in B(k)$ in a dense open subset of~$B$, the set $Z_b(k)$ is dense in $Z_b(\A_k)^{\Br(Z_b)}$.
\end{enumerate}
Then $Z(k)$ is dense in $Z(\A_k)^{\Br(Z)}$.
\end{thm}

The assumptions of Theorem~\ref{th:fibrationnf} imply that~$B$ is $k$\nobreakdash-rational.
Under the condition that~$B$ is $k$\nobreakdash-rational,
the first two assumptions of Theorem~\ref{th:fibrationnf} are expected to be superfluous (even under
weaker hypotheses on the generic fibre of~$p$ than rational connectedness,
see \cite[Corollary~9.23 (1)--(2)]{hwfibration}),
but removing
them altogether is a wide open problem, well connected with analytic number theory
(see \cite[\textsection9]{hwfibration}, \cite{hww}).
Removing~(2) while keeping~(1) might be within reach, though:

\begin{question}
In the statement of Theorem~\ref{th:fibrationnf}, can one dispense with the assumption
that~$p$ admit a rational section over~$\bar k$?
\end{question}

This would allow one
to replace ``supersolvable'' with ``solvable'' in the statement of Theorem~\ref{thm:arithhomspace}.
Indeed, in~\textsection\ref{subsec:geometry},
the cyclicity of the quotient $G_m/G_{m-1}$ plays a rôle  only to ensure
the existence of a rational
section of~$p$ over~$\bar k$ (see \cite[Proposition~3.3~(ii)]{hwzceh}).

\subsection{An application to Massey products}
\label{subsec:massey}

Theorem~\ref{thm:arithhomspace} has concrete applications, over number fields,
 beyond the inverse Galois problem:
for the homogeneous spaces that appear in its statement,
it turns the problem of deciding the existence of  a rational point
into the much more approachable question of deciding the non-vacuity
of the Brauer--Manin set.
In this way, Theorem~\ref{thm:arithhomspace} can be used to confirm, in the case of number fields,
 the conjecture of Mináč and Tân on the vanishing of Massey
products in Galois cohomology (see \cite{hwmassey}).
Indeed, this conjecture---which posits that 
for any field~$k$, any prime number~$p$,  any
integer $m\geq 3$ and any classes $a_1,\dots,a_m \in H^1(k,\Z/p\Z)$,
the  $m$\nobreakdash-fold
Massey product of $a_1,\dots,a_m$ vanishes if it is defined
(see \cite{minactanjems,minactantriple})---can be reinterpreted,
according to Pál and Schlank~\cite{PS16},
in terms of the existence of rational points on appropriate homogeneous spaces of~$\SL_n$ over~$k$
(with $n\gg 0$),
and it so happens that the geometric stabilisers of these homogeneous spaces
are finite and supersolvable.

\section{Rational curves on real algebraic varieties}
\label{sec:tight}

\subsection{A few questions}
Let~$X$ be a smooth variety over~$\R$.
The interplay between the topology of the $\ci$ manifold~$X(\R)$
and the geometry of the algebraic variety~$X$ lies at the core of classical real algebraic geometry.
One of the fundamental problems in this area consists in investigating which submanifolds of~$X(\R)$ can be approximated,
in the Euclidean topology, by Zariski closed submanifolds.
Even for
 $1$\nobreakdash-dimensional submanifolds, i.e.\ disjoint unions of~$\ci$ loops,
various phenomena---of a topological, Hodge-theoretic,
or yet more subtle
nature---can obstruct the existence of algebraic approximations
 (see \cite[\textsection4]{bwhodgereel1}).
In the case of $1$\nobreakdash-dimensional submanifolds,
however,
all known obstructions vanish
when~$X$ is rationally connected.
One can thus raise
the following questions,
in which $H_1^{\alg}(X(\R),\Z/2\Z)$ denotes the image of the cycle class map
$\CH_1(X) \to H_1(X(\R),\Z/2\Z)$ defined by
Borel and Haefliger \cite{borelhaefliger}.

\begin{questions}
\label{q:realrc}
Let~$X$ be a smooth, proper, rationally connected variety, over~$\R$.
\begin{enumerate}
\item Can all~$\ci$ loops in $X(\R)$ be approximated, in the Euclidean topology, by real loci of algebraic 
curves?
 or even
by real loci of rational algebraic
curves?
\item Is $H_1(X(\R),\Z/2\Z)=H_1^\alg(X(\R),\Z/2\Z)$?
 Is
 $H_1(X(\R),\Z/2\Z)$ generated by classes of rational algebraic curves on~$X$?
\end{enumerate}
\end{questions}

The first parts of Questions~\ref{q:realrc}~(1) and~(2)
are in fact equivalent to each other,
by the work of
Akbulut and King (see \cite[Theorem~6.8]{bwhodgereel2}), and were studied in a systematic fashion
in \cite{bwhodgereel1,bwhodgereel2}.
The second part of Question~\ref{q:realrc}~(1) is, however, as far as we know, genuinely stronger
than the second part of Question~\ref{q:realrc}~(2).
We note that in order to formulate the second part of
 Question~\ref{q:realrc}~(1) precisely,
it is better to work with
possibly non-injective~$\ci$ maps $\P^1(\R)\to X(\R)$
rather than with submanifolds of~$X(\R)$.
Indeed, there are examples of $\R$\nobreakdash-rational surfaces~$X$ and of~$\ci$ loops in~$X(\R)$
such that
the desired rational algebraic curves
necessarily have singular real points
(see \cite[Theorem~3]{kollarmangolteapprox}).

A specific motivation for Question~\ref{q:realrc}~(2) is its analogy with the following questions in complex geometry
raised by Voisin
\cite{voisinsomeaspects}
and by Kollár
\cite{kollarholomorphic}:

\begin{questions}
\label{q:voisinkollar}
Let~$X$ be a smooth, proper, rationally connected variety, over~$\C$.
Is the group $H_2(X(\C),\Z)$ generated
by homology classes
of algebraic curves? Is it generated by homology classes of rational algebraic curves?
\end{questions}

The two parts of Questions~\ref{q:voisinkollar} are in fact equivalent: Tian and Zong~\cite{tianzong} have shown
that the homology class of any algebraic curve on a rationally connected variety over~$\C$ is a linear
combination of homology classes of rational curves.  The real analogue of their result remains unknown
in general. Its validity is an interesting open problem.

The first parts of Questions~\ref{q:realrc}~(2) and of Questions~\ref{q:voisinkollar}
are in fact related by more than an analogy: if~$X$ is a smooth, proper, rationally connected variety
over~$\R$
such that
 $X(\R)\neq\emptyset$
and such that
 Questions~\ref{q:voisinkollar} admit a positive
answer for~$X_\C$,
then the equality
 $H_1(X(\R),\Z/2\Z)=H_1^\alg(X(\R),\Z/2\Z)$
is equivalent to the \emph{real integral Hodge conjecture} for $1$\nobreakdash-cycles on~$X$,
a property formulated and studied in \cite{bwhodgereel1,bwhodgereel2}.

In a different line of investigation around
the abundance of rational curves on rationally connected varieties,
many authors have considered the problem of finding rational curves
through a prescribed set of
points, or more generally through a curvilinear $0$\nobreakdash-dimensional subscheme,
on any smooth, proper, rationally connected variety~$X$.
Over the complex numbers,
such curves exist unconditionally (Kollár, Miyaoka, Mori, see \cite[Chapter~IV.3]{kollarbook}).
Over the real numbers, such curves exist under the necessary condition that
all the prescribed points that are real belong to the same connected component
of~$X(\R)$ (Kollár, see \cite{kollarlocal}, \cite{kollarspecialization}).
This problem can be generalised to one-parameter families:
 given a morphism $f:\sX \to B$ with rationally connected generic
fibre
between smooth and proper varieties,
where~$B$ is a curve,
one looks for sections of~$f$
whose restriction to
a given $0$\nobreakdash-dimensional subscheme of~$B$ is prescribed,
thus leading to
Questions~\ref{q:warc} below.
For simplicity of notation,
in the statement of
Question~\ref{q:warc}~(2),
this $0$\nobreakdash-dimensional subscheme of~$B$ 
is assumed to be reduced; there is however no loss of generality in doing this,
since jets of sections can be prescribed
at any higher order
by replacing~$\sX$ with a suitable
iterated blow-up (see \cite[Proposition~1.4]{Hassett}).

\begin{questions}
\label{q:warc}
Let~$B$ be a smooth, proper, connected curve over a field~$k_0$.
Let~$\sX$ be a smooth, proper variety over~$k_0$, endowed with
a flat morphism $f:\sX \to B$ with rationally connected generic
fibre.  Let $P \subset B$ be a reduced $0$\nobreakdash-dimensional subscheme.
Let $s:P \to \sX$ be a section of~$f$ over~$P$.
\begin{enumerate}
\item If $k_0=\C$, can~$s$ be extended to a section of~$f$?
\item If $k_0=\R$ and the map $s|_{P(\R)}: P(\R)\to \sX(\R)$ can be extended to a~$\ci$
 section of $f|_{\sX(\R)}:\sX(\R)\to B(\R)$, can then~$s$ be extended to a section of~$f$?
\end{enumerate}
\end{questions}
Let~$X$ be the generic fibre of~$f$ and~$k$ the function field of~$B$.
The existence of sections extending any given~$s$ as above is equivalent to the density of~$X(k)$
in the topological space $X(\A_k)=\prod_{b} X(k_b)$ of adelic points of~$X$, where
 the product runs over the closed points~$b$ of~$B$
and~$k_b$ denotes the completion of~$k$ at~$b$. This is the \emph{weak approximation} property.

The Graber--Harris--Starr theorem \cite{ghs} provides a positive answer to
Question~\ref{q:warc}~(1) when $P=\emptyset$ and it is a conjecture of
Hassett and Tschinkel that the answer to this question is in the
affirmative in general (see \cite{CTG},
\cite{HT}, \cite{Hassett}, \cite{tiancubic} for known results).
Particular cases of Question~\ref{q:warc}~(2) were first studied
by Colliot-Thélène~\cite{ctgroupes}, 
who conjectured the validity of weak approximation
(i.e.\ a positive answer to Question~\ref{q:warc}~(2)
even without assuming that
 $s|_{P(\R)}$ can be extended to a~$\ci$
 section of $f|_{\sX(\R)}$)
when~$X$ is birationally equivalent to a homogeneous space of a connected linear algebraic
group over~$k$, and proved his conjecture when the geometric stabilisers are trivial.
Scheiderer~\cite{Scheiderer}
then proved the same conjecture
 when the geometric stabilisers are connected.
Ducros~\cite{Ducros,DucrosCRAS} stated
Question~\ref{q:warc}~(2) in these exact terms, and gave a positive answer when~$X$ is a conic bundle surface,
or more generally when there exists a dominant map $X\to \P^1_k$
whose generic fibre is a Severi--Brauer variety.

\subsection{Tight approximation}
\label{subsec:tightapprox}

The main insight behind the proof of Theorem~\ref{thmintro:B} is the
observation that formulating a suitable common strengthening of
Questions~\ref{q:realrc} and Questions~\ref{q:warc}, through the notion of \emph{tight approximation},
 can render all of these
questions fully amenable to both the descent method and the fibration method.
We note that Questions~\ref{q:realrc} and Questions~\ref{q:warc} are somewhat orthogonal in spirit, insofar as
the former consider global constraints on curves lying on~$X$, while the latter
are aimed at local constraints.

The idea of establishing a descent method (resp.\ fibration method) for Question~\ref{q:warc}~(2)
already appeared in \cite{Ducros} (resp.~\cite{palszabo}), though in \emph{op.\ cit.}\ the implementations are subject
to miscellaneous restrictions.
The possibility of a descent method and a fibration method for
studying Questions~\ref{q:realrc}, however,
is new and turns out to require a shift in perspective  from single rationally connected
varieties to one-parameter families of such.

Let us illustrate how Questions~\ref{q:realrc}  need to be strengthened
for
a fibration argument to go through.  We start with a dominant morphism $p:X\to Y$ with rationally connected
generic fibre
between smooth, proper,
rationally connected
varieties, over~$\R$, and a $\ci$ loop $\gamma:\S^1 \to X(\R)$ that we want to approximate,
in the Euclidean topology, by a Zariski
closed submanifold of~$X$, assuming that we can solve the same problem on~$Y$ as well as on the fibres
of~$p$.  By assumption, we can approximate $p \circ \gamma : \S^1 \to Y(\R)$ by a~$\ci$
map $\xi:\S^1 \to Y(\R)$ with Zariski closed image.  The best we can hope to find, then, is
a $\ci$ loop $\tilde\gamma:\S^1 \to X(\R)$ arbitrarily close to~$\gamma$ and
such that $p \circ \tilde\gamma = \xi$.
We draw two conclusions:
\begin{enumerate}
\item If such a~$\tilde\gamma$ exists, the next and final step is not finding an algebraic approximation
for a $\ci$ loop in a fibre of~$p$, but, rather, considering
the algebraic curve~$B$ underlying~$\xi(\S^1)$,
viewing~$\tilde\gamma$ as a~$\ci$ section of the projection $(X \times_Y B)(\R) \to B(\R)$,
and looking for an algebraic section of $X \times_Y B \to B$ approximating~$\tilde\gamma$.
Thus, even when we start with just two real varieties~$X$ and~$Y$,
we need to consider one-parameter algebraic families of fibres of~$p$, rather than single fibres.
\item Consider the example where~$p$ is the blow-up of a surface~$Y$ at a real point~$b$
and~$\gamma$ meets $p^{-1}(b)(\R)$, transversally.
Then for any~$\tilde\gamma$ sufficiently close to~$\gamma$ in the Euclidean
topology, the loop $p \circ \tilde\gamma$ has to go through~$b$.
Hence~$\xi$ has to be required to go through~$b$
for a loop~$\tilde\gamma$ as above to exist.  Thus, a condition of weak approximation type must be considered
in conjunction with Questions~\ref{q:realrc} (as was already noted by Bochnak and Kucharz
 \cite{BKrat}).
\end{enumerate}
Let us now similarly contemplate a fibration argument in the context of Question~\ref{q:warc}~(2).
We assume
that $\sX \xrightarrow{f} B$ can be factored as $\sX\xrightarrow{p}\sY\xrightarrow{g} B$,
where the variety~$\sY$ is smooth and proper over~$\R$, the morphism~$p$ is dominant with rationally connected
generic fibre,
and~$g$ is flat.
Starting from a section $s:P\to \sX$ of~$f$ over~$P$ such that $s|_{P(\R)}$ can be extended to a~$\ci$
section~$s'$ of~$f|_{\sX(\R)}$, a positive answer to Question~\ref{q:warc}~(2) for~$g$ produces for us
a section~$\tau$ of~$g$ that extends $p\circ s$.
Let $\sZ = p^{-1}(\tau(B))$ and let $h:\sZ \to B$ denote the restriction of~$f$.
  At this point,
one would like to apply a positive answer to Question~\ref{q:warc}~(2) for~$h$
to obtain a section of~$h$ extending~$s$, thus completing the argument, as $\sZ \subseteq \sX$.
In order to do so, one needs to know that $s|_{P(\R)}:P(\R) \to \sZ(\R)$
can be extended to a~$\ci$ section of $h|_{\sZ(\R)}:\sZ(\R)\to B(\R)$.
However, the  map
 $h|_{\sZ(\R)}$
 in general even fails to be surjective.
To correct this problem, one should require, at the very least, that $\tau(B(\R))$ approximate,
 in the
Euclidean topology,
 the image of $p\circ s':B(\R)\to \sY(\R)$.
Thus, all in all,
an approximation condition in the Euclidean topology
has to be considered in conjunction with Question~\ref{q:warc}~(2).

The above discussion leads to the following definition.  (This definition slightly
differs from the one given
in \cite{bwtight}, which considers the more general question of approximating holomorphic maps by
algebraic ones, à la Runge, and which, as a consequence, is useful also for studying complex curves
on complex varieties, without reference to the reals;
however, all of the statements we make below are true with respect to either of the definitions.)

\begin{defn}
\label{def:tight}
Let~$B$ be a smooth, proper, connected curve over~$\R$.  A  variety~$X$ over $k=\R(B)$
satisfies the \emph{tight approximation} property if for any proper model $f:\sX\to B$ of~$X$ over~$B$
with~$\sX$ smooth over~$\R$,
 any reduced $0$\nobreakdash-dimensional subscheme $P\subset B$,
any section $s':P\to \sX$ of~$f$ over~$P$
and
 any~$\ci$ section $s:B(\R)\to \sX(\R)$ of $f|_{\sX(\R)}$
such that $s|_{P(\R)}=s'|_{P(\R)}$,
there exists a section $\sigma:B \to \sX$ of~$f$ such that
 $\sigma|_{P}=s'|_{P}$
and such that
 $\sigma|_{B(\R)}$ lies arbitrarily close
to~$s$ in the compact-open topology.
\end{defn}

Given a smooth, proper, rationally connected variety~$X$ over~$\R$,
the validity of the tight
approximation property for the variety obtained from~$X$ by extension of scalars from~$\R$ to~$\R(t)$
implies
positive answers to Questions~\ref{q:realrc} for~$X$.

The tight approximation property is (tautologically) a birational invariant, and it holds
for~$\P^n_k$
by a theorem of
Bochnak and Kucharz~\cite{BKrat}. (In \emph{op.\ cit.},
weak approximation conditions at complex points are ignored, but they create no additional difficulty.)
The next two results provide more examples of varieties satisfying tight approximation.

\subsection{Descent}

The following theorem
implements the descent method
for the tight approximation
property,
 in full generality (including non-abelian descent,
as formalised by Harari and Skorobogatov).
Its proof, given in \cite{bwtight},
 builds on the work of Scheiderer~\cite{Scheiderer} and, in the case where~$G$ is finite,
on an argument of Colliot-Thélène and Gille~\cite{CTG}.

\begin{thm}
\label{th:descenttight}
Let~$k$ be the function field of a real curve.
Let~$X$ be a smooth variety over~$k$.
Let~$G$ be a linear algebraic group over~$k$.
Let $f:Y\to X$ be a left torsor under~$G$.
Consider twists $f':Y'\to X$ of~$f$ by right torsors under~$G$, over~$k$.
If every such~$Y'$ satisfies the tight approximation property, then so does~$X$.
\end{thm}

\subsection{Fibration}

The next theorem implements the fibration method for the tight approximation property,
in full generality.
Its proof, contained in~\cite{bwtight}, makes essential use of
the weak toroidalisation theorem of Abramovich, Denef and Karu~\cite{ADK}
to establish a version of the Néron smoothening process
(as in \cite[3.1/3]{neronmodels})
 for higher-dimensional bases---the point being that in the discussion
at the beginning of~\textsection\ref{subsec:tightapprox}, the loop~$\tilde\gamma$ is easily seen
to exist once the morphism~$p$ is smooth along~$\gamma$ (see \cite[Lemma~6.11]{bwhodgereel2}).

\begin{thm}
\label{th:fibrationtight}
Let~$k$ be the function field of a real curve.
Let $p:Z\to B$ be a dominant morphism between smooth varieties over~$k$.
If~$B$ and the fibres of~$p$ above the rational points of a dense open subset of~$B$
satisfy the tight approximation property, then so does~$Z$.
\end{thm}

\subsection{Homogeneous spaces}

We are now in a position to sketch the proof of the following theorem,
which
in the ``constant case'', i.e.\ when the algebraic group
and the homogeneous space are both defined over~$\R$,
immediately implies
Theorem~\ref{thmintro:B}.

\begin{thm}
\label{th:homspacetight}
Homogeneous spaces of connected linear algebraic groups over the function field of a real curve
satisfy the tight approximation property.
\end{thm}

The proof of Theorem~\ref{th:homspacetight}
starts by noting that quasi-trivial tori over~$k$ are $k$\nobreakdash-rational,
hence satisfy the tight approximation property (since so does~$\P^n_k$).
Any torus~$T$ can be inserted into an exact sequence $1\to S\to Q\to T\to 1$ where~$S$ is a torus
and~$Q$ is a quasi-trivial torus.
As any twist of~$Q$ as a torsor remains isomorphic to~$Q$
(Hilbert's Theorem~90)
and hence satisfies the tight approximation property,
we deduce, by the descent method
(Theorem~\ref{th:descenttight}), that all
tori over~$k$ satisfy the tight approximation property.
Next, as every connected linear algebraic group over~$k$ is birationally equivalent to a relative torus
over a $k$\nobreakdash-rational variety (namely over the variety of maximal tori, when the algebraic
group is reductive),
we deduce, by the fibration method
(Theorem~\ref{th:fibrationtight}), that connected linear algebraic groups over~$k$ satisfy the
tight approximation property.  By descent
(Theorem~\ref{th:descenttight} again), it follows that homogeneous spaces
of connected linear algebraic groups over~$k$
satisfy the tight approximation property when they have a rational point.
Finally, it is a theorem of Scheiderer that homogeneous spaces of connected linear algebraic groups
over~$k$ satisfy the Hasse principle with respect to the real closures of~$k$,
so that if~$X$ denotes such a homogeneous space, then $X(k)\neq\emptyset$ whenever
a~$\ci$ section $s:B(\R)\to \sX(\R)$ as in Definition~\ref{def:tight} exists.
This completes the proof of Theorem~\ref{th:homspacetight}.

\subsection{Further comments}

Theorem~\ref{th:homspacetight} implies
that
 homogeneous spaces of connected linear algebraic groups over the function field of a real
curve satisfy
 weak approximation,
as conjectured by Colliot-Thélène.
Indeed, in
the notation of Definition~\ref{def:tight},
if~$X$ is such a homogeneous space
and~$P$ contains the locus of singular fibres of~$f$,
Scheiderer's work implies that $f^{-1}(b)(\R)$ is nonempty and connected
for all $b \in B(\R) \setminus P(\R)$,
so that a $\ci$ section $s:B(\R)\to\sX(\R)$ with
 $s|_{P(\R)}=s'|_{P(\R)}$ always exists.

The main open problem surrounding the notion of tight approximation is the following.

\begin{question}
\label{q:tightrc}
Let~$k$ be the function field of a real curve.
Do all rationally connected varieties over~$k$ satisfy the tight approximation
property?
\end{question}

Building on Theorem~\ref{th:descenttight}
and Theorem~\ref{th:fibrationtight},
the tight approximation property
is shown
in~\cite{bwtight}
to hold
for various classes of rationally connected varieties beyond homogeneous spaces of connected linear
algebraic groups.
For instance,
it holds for smooth cubic hypersurfaces of dimension~$\geq 2$
that are defined over~$\R$,
thus yielding, for such hypersurfaces, a positive answer
to (the second part of) Question~\ref{q:realrc}~(1).

Question~\ref{q:tightrc} is open for cubic surfaces over~$k$.
Even Question~\ref{q:warc}~(2) is open
when~$X$ is a cubic surface, although Question~\ref{q:warc}~(1) has an affirmative answer
in this case, by a theorem of Tian~\cite{tiancubic}.

In another direction,
Question~\ref{q:tightrc} is open for surfaces defined over~$\R$,
and so is (the second part of) Question~\ref{q:realrc}~(1).
By inspecting the birational classification of geometrically rational surfaces
and
using the fibration method (Theorem~\ref{th:fibrationtight}),
one can see that a positive answer to these questions
for surfaces defined over~$\R$
would follow from a positive answer for
del Pezzo surfaces of degree~$1$ or~$2$ defined over~$\R$.
In these cases,
it would suffice, by an application of the descent method (Theorem~\ref{th:descenttight}),
to know that for any real del Pezzo surface~$X$ of degree~$1$ or~$2$,
the universal torsors of~$X$, in the sense of Colliot-Thélène and
Sansuc~\cite{ctsandescent2}, are $\R$\nobreakdash-rational whenever they have a real point.
This last question, unfortunately, is very much open---even the unirationality of real
del Pezzo surfaces of degree~$1$ is unknown.  In fact, not a single example of a
minimal real del Pezzo surface
of degree~$1$ is known to be unirational.
For a description of these surfaces,
see \cite[\textsection5]{russo}.

Naturally, one hopes for the answer to Question~\ref{q:tightrc}
to be in the affirmative in general.
This conjecture would have a host of interesting consequences, among which:
a version of the Graber--Harris--Starr theorem over the reals (i.e.\ a positive answer to
Question~\ref{q:warc}~(2) when $P=\emptyset$);
Lang's widely open conjecture
from
 \cite{langrealplaces} that the function field of a real curve with no real point is~$C_1$
(see \cite[Corollary~1.5]{hogadixu} for the implication);
and the existence of a geometrically rational curve on any
 smooth, proper, rationally connected variety of dimension~$\geq 1$ over~$\R$.

This last consequence is a conjecture of Kollár, who showed the existence of rational curves on
those real rationally connected
varieties of dimension~$\geq 1$ that have real points
(see \cite[Remarks~20]{araujokollar}).  For real rationally connected
varieties with no real point, it is interesting to consider a weaker property:
the existence of a geometrically irreducible curve of even geometric
genus.  The latter can be reinterpreted
in terms of the
real integral Hodge conjecture
(see \cite{bwhodgereel1}).  Using Hodge theory and a real adaptation of Green's infinitesimal criterion for the density of Noether--Lefschetz loci, 
such curves of even genus can be shown
to exist on all real Fano threefolds (see \cite{bwhodgereel2}).
However, even
on smooth quartic hypersurfaces in~$\P^4_\R$, the existence
of geometrically rational curves remains a challenge, as well as
the mere existence of an absolute bound, independent of the chosen quartic hypersurface, on the minimal geometric genus of a geometrically irreducible curve of even geometric genus lying on such a hypersurface.

\section{Function fields of curves over \texorpdfstring{$p$\nobreakdash-adic}{𝑝-adic} fields}
\label{sec:otherfields}

\subsection{Some motivation: rational curves over number fields}

Even though the main questions about rational points of
rationally connected varieties
over number fields and over function fields of real curves
are still wide open,
the Brauer--Manin obstruction and the tight approximation property at least provide rather
satisfactory conjectural answers.
It would be highly desirable to obtain a similar conjectural picture for rational points
over other fields, for significant classes of
varieties---including, at a minimum, concrete criteria for the existence of rational points.

Over the field~$\Q(t)$,
this would encompass
questions about rational curves on rationally connected varieties over~$\Q$,
about which very little is known.  For example, it is unknown whether
 any rationally connected variety of dimension~$\geq 1$ over~$\Q$ that possesses
a rational point also contains a rational curve defined over~$\Q$.
Much more ambitiously, it is unknown whether
any such variety contains enough rational curves to imply the finiteness of the set of $R$\nobreakdash-equivalence
classes of rational points, a question asked in \cite[Question~10.12]{ctdegenerescences}.
(Known results on this problem are listed in \emph{loc.\ cit.})
As another example,
the regular inverse Galois problem over~$\Q$, which asks for the construction of
a regular Galois extension of~$\Q(t)$ with specified Galois group, and
which can be reinterpreted as a problem about
the existence of appropriate rational curves
on the homogeneous space~$\SL_n/G$ over~$\Q$,
is open even for finite nilpotent groups~$G$.
All of these problems are currently out of reach.

As a first step towards these questions, let us replace~$\Q$ with its completions and turn to rational points
over the field~$\Q_p(t)$ or over its finite extensions.

\subsection{Rational curves on varieties over \texorpdfstring{$p$\nobreakdash-adic}{𝑝-adic} fields}

In the constant case
 (that is, for varieties obtained by scalar extension from varieties defined over a $p$\nobreakdash-adic field,
i.e.\ a finite extension of~$\Qp$),
various existence results are known:
\begin{enumerate}
\item the regular inverse Galois problem over~$\Q_p$ has a positive
solution (Harbater \cite{harbatergalois}, by ``formal patching''; reproved and generalised
in different directions by Pop \cite{popgalois} and by Colliot-Thélène \cite{ctcovers}; see also
\cite{liuharbater},
\cite{moretbaillyrevetements},
\cite{kollarrcandpi1});
\item
for any smooth, proper, rationally connected variety~$X$ over a $p$\nobreakdash-adic field~$k$,
 Kollár \cite{kollarlocal,kollarspecialization}
has shown that the rational points of~$X$ fall into finitely many $R$\nobreakdash-equivalence classes,
and that there exist rational curves on~$X$, defined over~$k$, passing
 through any finite set of rational points of~$X$ that belong to the same $R$\nobreakdash-equivalence class
(with prescribed jets
of any given order at these points).
\end{enumerate}
This last statement concerns conditions of weak approximation type
that can be imposed on rational curves on rationally connected varieties over $p$\nobreakdash-adic fields.
It would be interesting to formulate
an
analogue, in this $p$\nobreakdash-adic context, of the surjectivity of the Borel--Haefliger
cycle class map
$\CH_1(X) \to H_1(X(\R),\Z/2\Z)$
(i.e.\ of Questions~\ref{q:realrc}~(2)).

We saw in~\textsection\ref{sec:tight}
that
in order to answer questions about homology classes of rational curves on real varieties, it can be useful
to consider more generally the tight approximation property,
for non-constant varieties over the function field
of a real curve.  By analogy, this gives incentive to investigate the possibility of
a $p$\nobreakdash-adic analogue of the tight approximation property for non-constant varieties
over the function field of a curve over a $p$\nobreakdash-adic field,
the validity of which would have consequences for
 a likely easier to formulate \emph{$p$\nobreakdash-adic integral Hodge conjecture} for $1$\nobreakdash-cycles
on varieties over $p$\nobreakdash-adic fields.

\subsection{Quadrics and other homogeneous spaces}
\label{subsubsec:quadrics}

In the non-constant case,
even the simplest varieties over~$\Q_p(t)$
lead to difficult problems when it comes to  their rational points.
For instance, it is only a relatively recent theorem of Parimala and Suresh~\cite{parimalasureshuinvariant1},
for $p\neq 2$, and of Leep~\cite{leep},
based on work of Heath-Brown~\cite{heathbrownbeforeleep},
for arbitrary~$p$,
that
every projective quadric of dimension~$\geq 7$ over $\Q_p(t)$ possesses a rational point.
(In the language of quadratic forms, ``the $u$\nobreakdash-invariant of $\Q_p(t)$ is equal to~$8$''.)
Many other articles have been devoted to local-global principles for varieties over function fields of
curves over $p$\nobreakdash-adic fields
(e.g.\ \cite{huyongweakapprox,hhorig,hhkuinvariant,ctpscmh,preetirost,hhkgalois,huyonghasse,pschar2,hararischeidererszamuely,hhktorsors,ctpssuperieur,harariszamuelylocglob,huyongcohohasse,ppscompos,hhkpsisrael,cthhkpszerocycles,hhkpcomparison,mehmetipatching,cthhkpstori,parimalasureshunitary,hkpirutka,tiany,cthhkpsreductive}).

A patching technique 
was developed by Harbater, Hartmann and Krashen
(``patching over fields'', a successor to formal patching),
and was applied to study rational points of homogeneous spaces over such fields.
It was used,
in~\cite{hhkuinvariant},
 to give another proof of the aforementioned theorem of Parimala and Suresh,
and, in \cite{ctpscmh}, to establish, more generally,
the local-global principle for the existence of rational points on
smooth projective quadrics of dimension~$\geq 1$ over~$\Q_p(t)$
(or over a finite extension of~$\Q_p(t)$),
 with respect
to all discrete valuations on this field,
when~$p$ is odd.

\subsection{Reciprocity obstructions}

Let~$k$ be a finite extension of~$\Q_p(t)$.  Let~$\Omega$ denote the set of equivalence classes of
discrete valuations (of rank~$1$) on~$k$ and, for~$v \in \Omega$, let~$k_v$ denote the completion of~$k$ at~$v$.
Let~$X$ be an irreducible, smooth and proper variety over~$k$.
We embed~$X(k)$ diagonally
into the
product topological space $\prod_{v \in \Omega} X(k_v)$,
which we shall also denote $X(\A_k)$ (recall that~$X$ is proper).

We now explain how,
 building on the work of Bloch--Ogus and of Kato, an analogue of the Brauer--Manin obstruction
can be set up in this context.
These ideas, which are due to Colliot-Thélène,
 appear in print, and are put to use,
in \cite[\textsection2.3]{ctpssuperieur},
in a very slightly different (equicharacteristic) situation.
We refer the reader to \emph{loc.\ cit.}\ for more details.
(The ``reciprocity obstructions'' of \cite[\textsection4]{hararischeidererszamuely}
are weaker than those we discuss here.)

Our goal is thus to define, in complete generality,
a closed subset $X(\A_k)^{\rec} \subseteq X(\A_k)$ containing~$X(k)$, using on the one hand
a reciprocity law coming from~$k$ and on the other hand an analogue of the Brauer group of~$X$.

Grothendieck's purity theorem for the Brauer group equates $\Br(X)$
with the unramified cohomology group $H^2_{\nr}(X/k,\Q/\Z(1))$.
We recall the definition of unramified cohomology:
 for any irreducible smooth variety~$V$ over a field~$K$ of characteristic~$0$
and any torsion Galois module~$M$ over~$K$,
the group $H^q_{\nr}(V/K,M)$
is 
the subgroup of the Galois cohomology group $H^q(K(V),M)$ consisting of those classes
whose residues along all codimension~$1$ points of~$V$ vanish.
  It is the unramified cohomology group $H^3_{\nr}(X/k,\Q/\Z(2))$ that will
serve as a substitute for $\Br(X)$ here.  (The shift in degree is explained by the fact that the field~$k$
 has cohomological
dimension~$3$ while number fields
have virtual cohomological dimension~$2$.)
For any field extension~$K/k$,
Bloch--Ogus theory
provides an evaluation map $H^3_{\nr}(X/k,\Q/\Z(2)) \to H^3(K,\Q/\Z(2))$,
$\alpha \mapsto \alpha(x)$
along any $K$\nobreakdash-point~$x$ of~$X$
(see  \cite{blochogus}).

Let~$\sB$ denote an irreducible normal proper scheme over~$\Z_p$ with function field~$k$.
In contrast with what happens over number fields, here it is not one reciprocity law that will play a rôle,
but infinitely many of them: one for each closed point of~$\sB$, for each such~$\sB$.
Namely, given any closed point $b \in \sB$, Kato \cite[\textsection1]{katohasse} has constructed a complex
\begin{align}
\label{eq:katocomplex}
H^3(k,\Q/\Z(2)) \longrightarrow \bigoplus_{\xi \in \sB_{1,b}} \Br(\kappa(\xi)) \longrightarrow \Q/\Z\rlap,
\end{align}
where~$\xi$ ranges
over
the set~$\sB_{1,b}$ of $1$-dimensional irreducible closed subsets of~$\sB$ that contain~$b$,
and where~$\kappa(\xi)$ denotes
the function
field of~$\xi$ (which is either a global field of characteristic~$p$ or a local field of characteristic~$0$).
The second arrow in~\eqref{eq:katocomplex} is the sum of the invariant maps from local class field theory
at the finitely many places of~$\kappa(\xi)$ that lie over~$b$.
The first arrow of~\eqref{eq:katocomplex} is induced by residue maps $\partial_v: H^3(k_v,\Q/\Z(2)) \to \Br(\kappa(\xi))$ constructed
by Kato in \emph{loc.\ cit.}, where~$v$ denotes the discrete valuation of~$k$ defined by~$\xi$.

For any $\alpha \in H^3_{\nr}(X/k,\Q/\Z(2))$, there are only finitely many $1$\nobreakdash-dimensional
irreducible closed subsets~$\xi$ of~$\sB$ such that
the map $X(k_v) \to \Br(\kappa(\xi))$, $x \mapsto \partial_v(\alpha(x))$ does not identically vanish,
if we denote by~$v$ the discrete valuation of~$k$ defined by~$\xi$
(see \cite[Proposition~2.7~(ii)]{ctpssuperieur} and note that for the proof given there,
it is enough to assume that a dense open subset of~$\sB$, rather than~$\sB$ itself, is a scheme over a field---an
assumption satisfied here).  As a consequence, it makes sense to define $X(\A_k)^{\rec}$ to be the set of
$(x_v)_{v \in \Omega} \in X(\A_k)$ such that
for any
 irreducible normal proper scheme~$\sB$ over~$\Z_p$ with function field~$k$,
 for any closed point $b \in \sB$, and for any $\alpha \in H^3_{\nr}(X/k,\Q/\Z(2))$,
the family 
$(\partial_v(\alpha(x_v)))_{\xi\in \sB_{1,b}} \in \bigoplus_{\xi \in \sB_{1,b}} \Br(\kappa(\xi))$
belongs to the kernel of the second arrow of~\eqref{eq:katocomplex}.
The fact that~\eqref{eq:katocomplex} is a complex immediately implies
that $X(k) \subseteq X(\A_k)^\rec$.

\subsection{Sufficiency of the reciprocity obstruction}

Although evidence is scarce,  the answer to the
following question might always be in the affirmative, as far as one knows:

\begin{question}
\label{q:bmonlyqpt}
Let~$k$ be a finite extension of~$\Q_p(t)$.
Let~$X$ be a smooth,
proper, rationally connected variety over~$k$.
If $X(\A_k)^{\rec}\neq\emptyset$,
does it follow that~$X(k)\neq\emptyset$?
\end{question}

Question~\ref{q:bmonlyqpt}
 has a positive answer
when~$X$ is a quadric and $p\neq 2$.
Indeed, we recall from~\textsection\ref{subsubsec:quadrics}
that
even $X(\A_k)\neq\emptyset$ then implies
 $X(k)\neq\emptyset$ (see \cite{ctpscmh}).
It
also has a positive answer when~$X$ is
 birationally equivalent to a torsor under a torus over~$k$.
This follows from the work of Harari, Scheiderer, Szamuely, Tian
\cite[Theorem~5.1]{harariszamuelylocglob},
\cite[\textsection0.3.1]{tianythese}
(modulo the comparison between the reciprocity obstruction defined
here and the reciprocity obstruction considered in these articles;
the latter is weaker but turns out to suffice to detect rational points on torsors under tori).
We note that there are examples of torsors under tori over~$k$ whose smooth compactifications~$X$ satisfy
$X(\A_k)^{\rec}=\emptyset$ while $X(\A_k)\neq\emptyset$
(see \cite[Remarque~5.10]{ctpssuperieur}).
Positive answers to Question~\ref{q:bmonlyqpt} are known
in various other cases in which~$X$ is birationally equivalent to a homogeneous space
of a connected linear algebraic group over~$k$.
For specific statements, we refer the reader
to the articles quoted in~\textsection\ref{subsubsec:quadrics}.
Question~\ref{q:bmonlyqpt} remains open in general for smooth compactifications of torsors under connected linear algebraic groups
over~$k$, for smooth compactifications of homogeneous spaces of~$\SL_n$ with finite
stabilisers, and for conic bundle surfaces over~$\P^1_k$.

Question~\ref{q:bmonlyqpt} focuses on the existence of rational points rather than on the density
of~$X(k)$ in~$X(\A_k)^{\rec}$ as the latter property is only known for
projective space (see \cite[Theorem~1]{artinwhaples})
and hence for varieties that
are rational as soon as they possess a rational point, such as quadrics.
For smooth compactifications of tori,
the density of~$X(k)$ in~$X(\A_k)^{\rec}$ is known to hold
off the set of discrete valuations of~$k$ whose residue field has
characteristic~$p$ (see~\cite[Theorem~5.2]{hararischeidererszamuely};
for the meaning of ``off'' here, see \cite[Definition~2.9]{wittenbergslc}).

To obtain more positive answers to Question~\ref{q:bmonlyqpt}, it is natural to wish
for flexible tools such as
general descent theorems and fibration theorems.
In the same way that 
introducing the tight approximation property
and
replacing
Question~\ref{q:warc}~(2)
with
Question~\ref{q:tightrc}
was a key step to obtain a problem that behaves well with respect to fibrations into rationally connected
varieties
(see the discussion in~\textsection\ref{subsec:tightapprox}),
it is likely that in order to obtain compatibility with descent and fibrations, one will have
to strengthen Question~\ref{q:bmonlyqpt} by incorporating into it a $p$\nobreakdash-adic analogue
of the  approximation condition
in the Euclidean topology
 that appears in Definition~\ref{def:tight}.
The main challenge, here, is to provide the correct formulation for such a $p$\nobreakdash-adic tight approximation property.

We note that in any case, a general fibration theorem has to lie deep, as it would presumably give a direct
route to the local-global principle for the existence of
rational points on smooth projective quadrics over~$k$ (so far unknown when $p=2$)
and hence to the computation of the $u$\nobreakdash-invariant of~$k$ (equal to~$8$;
see~\textsection\ref{subsubsec:quadrics}).
Indeed,
in the case of conics over~$k$,
this local-global principle follows from Tate--Lichtenbaum duality \cite{tatelichtenbaum};
applying a fibration theorem to
a general pencil of hyperplane sections of a fixed smooth projective quadric of dimension~$n \geq 2$
would allow one to deduce the general case by induction on~$n$.

\subsection{Further questions}

A good understanding of rational points of rationally connected varieties over function fields of curves
over $p$\nobreakdash-adic fields,
be it via
Question~\ref{q:bmonlyqpt}
or otherwise,
should shed light on concrete test questions such as the following:
\begin{questions}
\label{q:further}
Let~$p$ be a prime number and~$k$ be a finite extension of~$\Q_p(t)$.
\begin{enumerate}
\item Does
 the conjecture of Mináč and Tân on the vanishing of Massey
products in Galois cohomology hold for~$k$?
(See \textsection\ref{subsec:massey}
and
\cite{minactanjems,minactantriple}.)
\item Is there an algorithm that  takes as input a smooth, projective, rationally connected variety~$X$  over~$k$
and decides whether~$X$ has a rational point?
\end{enumerate}
\end{questions}
One might approach the first of these questions
by trying to mimic~\cite{hwmassey} over~$k$,
which would require making progress
on the arithmetic, over~$k$, of homogeneous spaces
of~$\SL_n$ with finite supersolvable geometric stabilisers.

To put the second question in perspective, let us recall what is known about
 algorithms for deciding the existence of rational points
on arbitrary varieties (``Hilbert's tenth problem'') over various fields of
interest.  Over~$\Q$ or~$\C(t)$, the existence of such an algorithm
is an outstanding open problem.  Denef~\cite{denefreal}
showed that
 over~$\R(t)$, such an algorithm does not exist.
His method was extended to prove that there is no such algorithm
over~$\Q_p(t)$ (Kim and Roush \cite{kimroush},
completed by Degroote and Demeyer \cite{degrootedemeyer}), over any finite extension of~$\R(t)$
that possesses a real place
(Moret-Bailly \cite{moretbaillyhilbertreal}),
or, when $p\neq 2$,  over any finite extension of~$\Q_p(t)$
(Eisenträger \cite{eisentragerpadic},
Moret-Bailly \cite{moretbaillyhilbertreal}).
In addition, over number fields,
it is known that
 restricting from arbitrary varieties to smooth projective varieties makes no difference
(see \cite[\textsection II.7]{smorynski}, \cite[Theorem~1.1~(i)]{poonenexistence}).
Restricting to smooth, projective, rationally connected varieties, however, does make a drastic difference:
Question~\ref{q:further}~(2) might well have an affirmative
answer for all of the fields just mentioned.
Over~$\C(t)$, this is  trivially so, by the Graber--Harris--Starr theorem.
Over~$\R(t)$,
a positive answer
to Question~\ref{q:further}~(2)
 would follow from 
a positive answer to Question~\ref{q:tightrc}.
Indeed, in the notation of Definition~\ref{def:tight},
if~$X$ satisfies the tight approximation property, then~$X$ has a rational point if and only if $f|_{\sX(\R)}$
admits a~$\ci$ section, a property that can be decided algorithmically.
Over number fields, as was observed by Poonen \cite[Remark~5.3]{poonenheuristics},
a positive answer
to Question~\ref{q:further}~(2) would
 follow from the conjecture that rational points are always
dense in the Brauer--Manin set.
It seems likely that a positive
answer to Question~\ref{q:bmonlyqpt}
would similarly imply a positive answer to Question~\ref{q:further}~(2).
To mimic Poonen's argument,
one runs into the difficulty that the elements of $H^3_{\nr}(X/k,\Q/\Z(2))$ are harder to describe
than those of $H^2_{\nr}(X/k,\Q/\Z(1))=\Br(X)$, whose interpretation  in terms of Azumaya algebras
is a key point in \emph{loc.\ cit.}; however, this can be remedied by viewing $H^3_{\nr}(X/k,\Q/\Z(2))$,
using Bloch--Ogus theory,
as the group of global sections of the Zariski sheaf associated with the presheaf
$U \mapsto H^3_{\et}(U,\Q/\Z(2))$, and describing $H^3_{\et}(U,\Q/\Z(2))$ via \v{C}ech cohomology.

\subsection{Other fields}
There are a number of other fields over which a better understanding of rational points of rationally
connected varieties would be valuable.  One of the simplest example is the fraction field $k=\C((x,y))$
of the ring of formal power series~$\C[[x,y]]$, which can be seen as a first step before considering
function fields of complex surfaces.  This field
presents both local and global features, and a reciprocity
obstruction
can again be defined (in terms of the unramified Brauer group---recall that~$k$ has cohomological
dimension~$2$).  This obstruction
was used in \cite{ctpssuperieur} to produce the first example of a torsor~$Y$ under
a torus, over~$k$, such that $Y(k)=\emptyset$ but $Y(k_v)\neq\emptyset$ for every
discrete valuation~$v$ on~$k$.
The analogues of Question~\ref{q:bmonlyqpt} and of Questions~\ref{q:further} can be asked
over this field too.  It is not known, however, whether
the reciprocity obstruction explains the absence of rational points on smooth  proper varieties
that are birationally equivalent to torsors
under tori over~$k$ (though see \cite[Corollaire~4.4]{izquierdodualitecxy}
for a closely related result involving possibly ramified Brauer classes).
We refer the interested reader to
\cite{ctop,ctgp,izquierdodualitecxy,izquierdolageometric}
for the state of the art.

\bibliographystyle{amsalpha}
\bibliography{icm}
\end{document}